\documentclass[11pt,english]{article}
\usepackage[T1]{fontenc}
\usepackage[latin1]{inputenc}
\usepackage{amssymb}

\makeatletter

\newcommand{\binom}[2]{{#1 \choose #2}}

\title{The Semigroup of a Word}
\author{Peter M. Higgins \& Norman R. Reilly}
\date{}
\def\qed{\quad\vrule height4.17pt width4.17pt depth0pt}

\makeatother

\usepackage{babel}
\begin{document}

\title{Embedding in a finite 2-generator semigroup}

\author{Peter M. Higgins, University of Essex U.K.}
\maketitle
\begin{abstract}
We augment the body of existing results on embedding finite semigroups
of a certain type into 2-generator finite semigroups of the same type.
The approach adopted applies to finite semigroups the idempotents
of which form a band and also to finite orthodox semigroups. 
\end{abstract}

\section{Introduction}

In this paper we will be concerned with the possibility of embedding
a finite semigroup $S$ into a finite $2$-generated semigroup $T$
that shares properties with $S$. In particular we show that any finite
orthodox semigroup $S$ may be embedded in a finite orthodox semigroup
$T$ generated by two group elements and that any finite orthodox
monoid $S^{1}$ may be embedded as a semigroup into a finite 2-generated
orthodox monoid $T$ whose subband of idempotents satisfies the same
semigroup identities. Prior to that we prove that if $S^{1}$ is a
finite monoid whose idempotents $E(S^{1})$ form a subsemigroup, then
$S^{1}$ may be embedded in a 2-generated finite monoid $T$ whose
idempotents also form a subsemigroup and belong to the same variety
of bands. For background on semigroups we refer to standard texts
such as {[}4{]} or {[}5{]}. 

Any semigroup $S$ may be embedded in the full transformation semigroup
$T=T_{S^{1}}$ (we shall sometimes write $S\leq T$ to denote that
$S$ is a subsemigroup of $T$). Since this natural `Cayley' embedding
preserves finiteness, it follows at once that any finite semigroup
$S$ embeds in the (regular) 3-generator semigroup $T_{n}$, where
$n=|S^{1}|$. We denote the corresponding semigroups of partial transformations
on a set $X$ by $PT_{X}$ and if $|X|=n$ we write this as $PT_{n}$.

In 1952 Trevor Evans proved in {[}2{]} that any countable semigroup
embeds in a 2-generator semigroup although that fact is implicit in
the paper {[}11{]} of Sierpinski published (in French) in 1935 where
it was shown that any countably infinite collection of mappings in
$T_{X}$ embeds in a 2-generator subsemigroup of $T_{X}$. The first
explicit proof that a finite semigroup may be embedded in a 2-generated
finite semigroup dates from 1960 and is due to B.H. Neumann {[}10{]}
who employed a wreath product construction. The short proof of this
fact recorded here however is indicative of the approach of the present
paper.

\textbf{Theorem 1.1} Any finite semigroup $S$ may be embedded in
a finite semigroup $T=\langle\alpha,\beta\rangle$ where $\alpha$
is an idempotent and $\beta$ is a nilpotent.

\emph{Proof }Without loss we assume that $S=S^{1}=\{\alpha_{0},\alpha_{1},\cdots,\alpha_{n-1}\}$
with $S\leq T_{X}$ for some finite set $X$ and where we take $\alpha_{0}=\iota$,
the identity mapping, in this instance with domain $X$. Our semigroup
$T\leq PT_{Z}$ where $Z=X\times\{0,1,2,\cdots,n\}$. We also put
$\alpha_{n}=\iota$. The designated generators $\alpha$ and $\beta$
are defined as follows: 
\[
(x,i)\cdot\alpha=(x\cdot\alpha_{i},0)\,(0\leq i\leq n)
\]
\[
(x,i)\cdot\beta=(x,i+1)\,(0\leq i\leq n-1).
\]
 In particular $\beta^{n+1}=0$, the empty mapping and $\alpha$ is
idempotent: 
\[
(x,i)\cdot\alpha^{2}=(x\cdot\alpha_{i},0)\cdot\alpha=(x\cdot\alpha_{i}\alpha_{0},0)=(x\cdot\alpha_{i},0)=(x,i)\cdot\alpha.
\]
 Hence $T$ is generated by an idempotent $\alpha$ together with
a nilpotent $\beta$. Now put $\lambda=\beta^{n}\alpha\in T$. Then
dom $\lambda=X\times\{0\}$ and 
\[
(x,0)\cdot\lambda=(x,0)\cdot\beta^{n}\alpha=(x,n)\cdot\alpha=(x\cdot\alpha_{n},0)=(x,0)
\]
 so that $\lambda=\iota|_{X\times\{0\}}$. Put $\gamma_{i}=\lambda\beta^{i}\alpha$
$(0\leq i\leq n-1)$; then dom $\gamma_{i}=X\times\{0\}$ and 
\[
(x,0)\cdot\gamma_{i}=(x,0)\cdot\lambda\beta^{i}\alpha=(x,0)\cdot\beta^{i}\alpha=(x,i)\cdot\alpha=(x\cdot\alpha_{i},0).
\]
 It follows that the mapping where $\alpha_{i}\mapsto\gamma_{i}$
is a monomorphism of $S$ into $T$, as required. $\qed$

It is not possible however to embed an arbitrary finite semigroup
into a finite semigroup generated by two idempotents as it is easy
to prove that any semigroup (finite or not) generated by two idempotents
has at most six idempotents and also does not contain a three-element
chain. A complete description of semigroups generated by two idempotents
has been provided by Benzaken and Mayr {[}1{]}.

In {[}7{]} Margolis showed that a finite semigroup $S$ may be embedded
in a 2-generated semigroup $T$ that is a Rees matrix semigroup $M(S)$
over $S$ with a cyclic group adjoined as group of units. This allowed
the conclusion that if all the subgroups of $S$ were abelian (nilpotent,
solvable, etc.), then you can embed $S$ into a 2-generator semigroup
$T$ with $T$ satisfying the same restriction on subgroups as $S$.
The construction idea was used in {[}6{]} to show that a compact metric
semigroup may be embedded in a 2-generator compact monoid. Moreover
it is implicit in {[}7{]} that any (finite) $n$-generated semigroup
$S$ may be embedded in a (finite) semigroup $T$ generated by $n+1$
idempotents, from which it follows that any finite semigroup $S$
may be embedded in a finite semigroup generated by three idempotents. 

Although not the principle result in their paper, in {[}8{]} McAlister,
Stephen and Vernitski obtained a direct embedding of $T_{n}$ into
a 2-generator subsemigroup of $T_{n+1}$. Although they then move
on to the question of inverse semigroups (discussed below), their
construction implies the following result.

\textbf{Theorem 1.3 }Any finite semigroup may be embedded in a 2-generated
semigroup that is finite and regular.

It is enough to prove the result for $T_{n}$$(n\geq3)$ and in {[}8{]}
McAlister et. al. embed $T_{n}$ in a semigroup $S=\langle\alpha,\beta\rangle\leq T=T_{n+1}$.
We write the idempotent of defect $1$ in which $i\mapsto j$ $(i\neq j)$
as $\binom{i}{j}$. Using this notation, the generator $\beta$ is
the $(n+1)$-cycle $\beta=(1\,2\,\cdots\,n\,n+1)$ while $\alpha=(1\,2)\binom{n}{n+1}$,
a product of a transposition and an idempotent of defect $1$. That
$S$ contains a copy of $T_{n}$ then follows from a series of easily
verified facts: 

$\bullet$ The map $\varepsilon=\alpha^{2}=\binom{n}{n+1}$ is an
idempotent of defect $1$;

$\bullet$ for any $\gamma\in\varepsilon T\varepsilon$, consider
the restriction $\gamma|_{\{1,2,\cdots,n-1,n+1\}}$: this defines
an isomorphism of $\varepsilon T\varepsilon$ onto $T_{n}$ with base
set $\{1,2,\cdots,n-1,n+1\}$;

$\bullet$ $T_{n}$ is generated by the set consisting of the $n$-cycle
$(1\,2\,\cdots n-1\,n+1)$, the transposition $(1\,2)$ and the idempotent
of defect $1$, $\binom{n-1}{n+1}$;

$\bullet$ taking inverse images of these three mappings under the
isomorphism results in a set of three generators of $\varepsilon T\varepsilon$,
which are respectively $\kappa=(1\,2\,\cdots n-1\,n+1)\binom{n}{1}$,
$\alpha$, and the idempotent of defect $2$, $\phi=\binom{n}{n+1}\binom{n-1}{n+1}$.

$\bullet$ finally we note that $\varepsilon=\alpha^{2},\kappa=\varepsilon\beta\varepsilon$,
and $\phi=\beta\varepsilon\beta^{-1}\varepsilon$, and so $T_{n}\cong\varepsilon T\varepsilon\leq S$.

This concludes the proof in {[}8{]} that any finite semigroup may
be embedded in a finite semigroup that is generated by a pair of group
elements. (Note there are two minor corrections: the paper says that
$(n-1)\cdot\kappa=n$ when it should say that $(n-1)\cdot\kappa=n+1$
and $\beta$ is listed as one of the three generators of $\varepsilon T\varepsilon$
when it should say $\alpha$.)

\emph{Proof of Theorem 1.3 }To complete the proof we need only observe
that the semigroup $S$ is indeed regular. First note that 
\[
\varepsilon T\varepsilon\leq S\Rightarrow\varepsilon^{2}T\varepsilon^{2}\subseteq\varepsilon S\varepsilon\Rightarrow\varepsilon T\varepsilon\subseteq\varepsilon S\varepsilon\subseteq\varepsilon T\varepsilon,
\]
 so that 
\[
\varepsilon T\varepsilon=\varepsilon S\varepsilon=\alpha(\alpha S\alpha)\alpha\subseteq\alpha S\alpha=\alpha^{3}S\alpha^{3}\subseteq\varepsilon S\varepsilon
\]
 giving equality throughout and in particular that $\alpha S\alpha\cong T_{n}$
is a regular subsemigroup of $S$. 

Now take any $\gamma\in S$. Either $\gamma\in\langle\beta\rangle,$
and so $\gamma$ is a (regular) group element or, since $\alpha=\alpha^{3}$,
we may write $\gamma=\beta^{t}\sigma\beta^{s}$ for some $\sigma\in\alpha S\alpha$
and $0\leq t,s\leq n$. Taking any inverse $\sigma'\in V(\sigma)$
we may now check that $\beta^{-s}\sigma'\beta^{-t}\in V(\gamma)$.
Therefore the semigroup $S$ is indeed regular. $\qed$

Equally, the construction in {[}7{]} also preserves regularity and
so Theorem 1.3 is also implicit in the Margolis paper. In {[}3, Theorem
4.1{]}, Hall gives a result of C.J. Ash, which shows that any countable
inverse semigroup may be embedded in an inverse semigroup with two
generators and any finite inverse semigroup may be embedded in a finite
inverse semigroup that is generated as an inverse semigroup by two
generators. (In {[}8{]} it is shown that any finite inverse semigroup
may be embedded in a finite inverse semigroup that is generated as
a semigroup by two generators.) The construction we introduce here
is inspired by the model of Ash. We have one principal generator that
contains copies of all the mappings in $S$, the semigroup to be embedded,
while the second generator is a cycle. The domain and range of the
principal generator then consists of many copies of the base interval,
which are distributed among the cycle of intervals in such a way that
unwanted products, which might spoil the embedding, are avoided in
the mappings that are to be simulated.

\section{Mian-Chowla property }

The base set of the $2$-generator transformation semigroup $T$ will
consist of a cycle of a large number of copies of the underlying interval
on which act the members of the semigroup $S$, which is to be embedded
in $T$. However, the action of our principal mapping $\alpha$ that
simulates all the members of $S$ will be confined to a relatively
small number of sparsely spaced intervals. This will ensure that unwanted
products do not arise in the construction.

To this end, let $S=\{\alpha_{1},\cdots,\alpha_{n}\}$ be a finite
semigroup with $S$ defined by partial transformations on a finite
base set $X$. Since we are interested in embedding $S$ into a\emph{
$2$}-generator semigroup $T$ sharing some of the same properties
as $S$, we may assume that $n\geq3$. Moreover, without loss we may
assume that $S$ does not contain the empty mapping.

In order to make our construction free of unwanted non-zero products,
we make use of the following sequence of numbers, first introduced
in {[}9{]}.

\textbf{Definition 2.1 }The Mian-Chowla (MC) sequence is the sequence
of non-negative integers $m_{0},m_{1},\cdots$ recursively defined
as follows. Set $m_{0}=0$; for $i\geq1$ define $m_{i}$ to be the
least integer exceeding $m_{i-1}$ such that each difference between
distinct integers in the sequence $m_{0},m_{1},\cdots,m_{i}$ is unique.

\textbf{Remarks 2.2 }The recursive step of the MC sequence is well-defined
as by choosing a sufficiently large integer we may find some $m$
such that each difference $m-m_{j}$ has not appeared previously among
the differences of pairs taken from the sequence: indeed it is clear
that $m_{i}\leq2m_{i-1}+1$ so that $m_{i}\leq2^{i}-1$. The MC sequence
begins: 
\[
0,1,3,7,12,20,30,44,65,80,96,122,147,181,203,251,289,\cdots
\]
The recursive rule of definition of the MC sequence is often formulated
in the equivalent form that $m_{i}$ is the least integer such that
the list of all pairwise sums, $m_{j}+m_{k}$ for $j,k\leq i$, has
no repeats. Note that under this alternative formulation, $j=k$ is
not forbidden. 

In Section 3 we shall work with this particular sequence in our construction:
$m_{i}$ will denote the member of the MC sequence indexed by $i$.
However, the results will apply to \emph{any }strictly increasing
sequence of integers with the MC property, meaning that no number
appears as a difference between distinct members more than once. There
are of course any number of such sequences: for example the sequence
$k^{n}$, $n=0,1,2,\cdots$, for any base $k\geq2$ possesses the
MC property. Moreover the MC property is inherited by subsequences.
In Section 4 we shall also call upon the following specific fact.

\textbf{Lemma 2.3 }For $i,j,k,l\leq n$, if $i\geq j$ and $k\geq l$
then $(2^{i}+2^{j})-(2^{k}+2^{l})=2^{n}+2^{0}$ implies that $i=n,$
$l=0$ and $j=k$. 

\emph{Proof }If $i\leq n-1$ then $2^{i}+2^{j}\leq2\cdot2^{n-1}=2^{n}$
and the equation cannot hold. Hence $i=n$, giving $2^{j}-(2^{k}+2^{l})=1$.
Hence $j\geq2$ and since both sides of the equation are odd, it follows
that $l=0$, and so $j=k$. $\qed$

\textbf{Remark 2.4 }Unfortunately, the MC sequence lacks the corresponding
property as for example: 
\[
44+65=109=96+12+1+0\Leftrightarrow m_{7}+m_{8}=m_{10}+m_{4}+m_{1}+m_{0}
\]
\[
\Leftrightarrow(m_{7}+m_{8})-(m_{4}+m_{1})=m_{10}+m_{0}.
\]
Suppose that $M=m_{0},m_{1},\cdots,m_{n}$ is a (strictly increasing)
MC sequence of non-negative integers and put $m=1+m_{n}$. For any
set $A\subseteq M$ and $r\in\mathbb{Z}$, let us write $A+r=\{(a+r)\,\mbox{(mod \ensuremath{m}),\,\ensuremath{a\in A}\}}$.
Suppose that $|A|\geq3$ and $A+r\subseteq M$ with $r\not\equiv0$
(mod $m$). Without loss we may assume that $1\leq r\leq m-1$. By
hypothesis, for each $m_{i}\in A$, $(m_{i}+r)$ (mod $m$) $=m_{j}$
for some $0\leq j\leq n$. It follows that either $m_{j}-m_{i}=r$
or if $(m_{i}+r)$ (mod $m$$)=m_{i}+r-m,$ then $m_{j}-m_{i}=r-m$.
Let $m_{a},m_{b},m_{c}$ be three pairwise distinct members of $A$.
Consider, modulo $m,$ each of $m_{a}+r,m_{b}+r$ and $m_{c}+r$.
It now follows that for at least two of $m_{a},m_{b},m_{c}$, let
us say $m_{a}$ and $m_{b}$, there exist $m_{j},m_{k}\in M$ such
that $m_{j}-m_{a}=m_{k}-m_{b}$, contrary to the MC condition. Hence
we conclude:

\textbf{Lemma 2.5. }Let $M=m_{0},m_{1},\cdots,m_{n}$ be a finite
strictly increasing sequence of non-negative integers with the MC
property and put $m=1+m_{n}$. Suppose that $A\subseteq M$ is such
that $(A+r)$ (mod $m$) $\subseteq M$ for some $r\not\equiv0$ (mod
$m$). Then $|A|\leq2$. 

\section{Embedding in a semigroup generated by a nilpotent and a cycle}

In this section we construct a general embedding of a finite semigroup
$S$ into a 2-generated finite semigroup $T$, which preserves the
property that the idempotents form a subsemigroup. 

We will make use here of the easily proved result that in the presence
of the band identity $x=x^{2}$, any \emph{heterotypical identity
}$\phi$ (one in which a variable appears on one side only) implies
the identity $x=xyx$. It follows that any band satisfying $\phi$
is a rectangular band.

Let $S$ be a finite semigroup $S=\{\alpha_{1},\alpha_{2},\cdots,\alpha_{n}\}$.
We shall take $S$ to be a subsemigroup of $PT_{X}$, where $X$ is
a finite base set. We may also assume that the domain of each $\alpha_{i}$
is not empty. In the following construction we could replace the set
of mappings $\{\alpha_{i}\}$ by any generating set of $S$ but for
simplicity of notation we work with $S$ as the generating set for
$S$.

Let $\{m_{i}\}_{i\geq0}$ denote the MC sequence and let $Z=X\times\{0,1,2,\cdots,m_{2n-1}\}$.
Taking addition modulo $m=1+m_{2n-1}$, we take one generator of our
containing $2$-generator semigroup $T$ to be $\beta$ where:
\begin{equation}
(x,i)\cdot\beta=(x,i+1)\,\,(0\leq i\leq m_{2n-1})
\end{equation}
Since $\beta$ is a cycle, the notation $\beta^{r}$ is meaningful
for all integers $r$. We next specify the domain and range of our
second generator $\alpha$: dom $\alpha$ is contained in the union
of the $n$ copies of $X$, $Y_{i}=X\times\{m_{i}\}$ $(n\leq i\leq2n-1)$
while the range $Z\alpha$ is a subset of a second union of another
$n$ copies of $X$: $X\times\{m_{i}\}$ $(0\leq i\leq n-1)$. We
define the action of $\alpha$ on the \emph{interval} $Y_{n+j}=X\times\{m_{n+j}\}$
as we shall call it as: 
\begin{equation}
(x,m_{n+j})\cdot\alpha=(x\cdot\alpha_{j},m_{j})\,\,\,(0\leq j\leq n-1)
\end{equation}

\textbf{Definition 3.1 }Let $T=\langle\alpha,\beta\rangle$, with
$\alpha,\beta$ defined as in (1) and (2).

\textbf{Lemma 3.2 }The generators $\alpha$ and $\beta$ of $T$ satisfy
$\alpha^{2}=0$ and $\beta^{m}=\iota$, where $m=1+m_{2n-1}$. For
each $\gamma\in T$ and $0\leq i\leq n-1$ there exists some $0\leq j\leq n-1$
such that $(X\times\{i\})\gamma\subseteq X\times\{j\}$; moreover
if $(x,i)\cdot\gamma,\,(x',i')\cdot\gamma\in X\times\{j\}$ then $i=i'$.

\emph{Proof} The first two facts follow respectively from (2) for
$\alpha$ and from (1) for $\beta$. The claims in the second sentence
follow for $\gamma=\alpha,\beta$ as each mapping is one-to-one on
second components whence, by induction on the length of the product,
the same follows for an arbitrary product $\gamma$ of these two generators.
$\qed$

\textbf{Lemma 3.3 }Let $\gamma\leq_{{\cal J}}\alpha\beta^{r}\alpha$.
Then dom $\gamma\subseteq X\times\{i\}$ for some $i$ such that $0\leq i\leq m-1$.

\emph{Proof }First suppose that $\gamma=\rho\lambda\sigma$ with dom
$\lambda\subseteq X\times\{j\}$ say and that $(x,i)\in$ dom $\rho\lambda$
so that $(x,i)\cdot\rho\in X\times\{j\}$. It follows from Lemma 3.2
applied to $\rho$ that dom $\rho\lambda\subseteq X\times\{i\}$ and
then since dom $\rho\lambda\sigma\subseteq$ dom $\rho\lambda$, we
obtain dom $\gamma\subseteq X\times\{i\}$. Therefore it is enough
to prove the claim for a mapping $\gamma$ of the form $\gamma=\alpha\beta^{r}\alpha$.
Since dom $\gamma\subseteq$ dom $\alpha$, it follows that each member
of dom $\gamma$ has the form $(x,m_{n+j})$ for some $0\leq j\leq n-1$.
We then obtain: 
\begin{equation}
(x,m_{n+j})\cdot\alpha\beta^{r}\alpha=(x\cdot\alpha_{j},m_{j})\cdot\beta^{r}\alpha=(x\cdot\alpha_{j},(m_{j}+r)\,(\mbox{mod \ensuremath{m)}})\cdot\alpha
\end{equation}
Again by definition of $\alpha$ we infer that $m_{j}+r\equiv m_{n+l}$
(mod $m$) for some $0\leq l\leq n-1$. Now suppose that $(x',m_{n+j'})\in$
dom $\gamma$; by (3) we deduce that $m_{j'}+r\equiv m_{n+l'}$ (mod
$m$) say, so that $m_{n+l'}-m_{j'}\equiv m_{n+l}-m_{j}\equiv r$
(mod $m$). Since $0\leq m_{j},m_{j'}<m_{n+l},m_{n+l'}\leq m-1$,
it follows that these congruences imply the corresponding equalities
and that $r\not\equiv0$ (mod $m$). By the MC property however we
conclude that $j=j'$ and $l=l'$. In particular, dom $\gamma\subseteq X\times\{m_{n+j}\}$,
giving the required conclusion. $\qed$

\textbf{Lemma 3.4 }Define the mapping $\lambda_{0}=(\beta^{m_{n}}\alpha)^{2}$.
Then\textbf{ }$\lambda_{0}=\iota|_{X\times\{0\}}$. 

\emph{Proof} From the definition of $\lambda_{0}$ we obtain 
\[
(x,0)\cdot\lambda_{0}=(x,0)\cdot(\beta^{m_{n}}\alpha)^{2}=(x,m_{n})\cdot\alpha\beta^{m_{n}}\alpha=(x\cdot\alpha_{0},m_{0})\cdot\beta^{m_{n}}\alpha=
\]
\[
(x,0)\cdot\beta^{m_{n}}\alpha=\cdots=(x,0).
\]
The result now follows from this together with Lemma 3.3. $\qed$

\textbf{Lemma 3.5 }The semigroup $T=\langle\alpha,\beta\rangle$ contains
each of the mappings $\lambda_{i,j,k}=\lambda(\alpha_{i},j,k)$ where
dom $\lambda_{i,j,k}\subseteq X\times\{j\},$ ran $\lambda_{i,j,k}\subseteq X\times\{k\}$
and $(x,j)\cdot\lambda_{i,j,k}=(x\cdot\alpha_{i},k)$ $(0\leq i\leq n-1,\,0\leq j,k\leq m-1)$.

\emph{Proof }We verify that $\lambda(\alpha_{i},j,k)=\beta^{-j}\lambda_{0}\beta^{m_{n+i}}\alpha\beta^{k-m_{i}}$.
Consider $(x,t)$ with $t\not\equiv j$ (mod $m$). Then $(x,t)\cdot\beta^{-j}=(x,t-j)\not\in X\times\{0\}$
so that by Lemma 3.4, $(x,t-j)\not\in$ dom $\lambda_{0}$. It follows
that dom $\beta^{-j}\lambda_{0}\beta^{m_{n+i}}\alpha\beta^{k-m_{i}}\subseteq X\times\{j\}$.
Next take $(x,j)\in X\times\{j\}$: 
\[
(x,j)\cdot\beta^{-j}\lambda_{0}\beta^{m_{n+i}}\alpha\beta^{k-m_{i}}=(x,0)\cdot\lambda_{0}\beta^{m_{n+i}}\alpha\beta^{k-m_{i}}=(x,0)\cdot\beta^{m_{n+i}}\alpha\beta^{k-m_{i}}=
\]
\[
=(x,m_{n+i})\cdot\alpha\beta^{k-m_{i}}=(x\cdot\alpha_{i},m_{i})\cdot\beta^{k-m_{i}}=(x\cdot\alpha_{i},k).
\]
 Therefore $\lambda(\alpha_{i},j,k)\in T$. $\qed$

\textbf{Theorem 3.6 }(Structure of $T$) 

(i) The monoid $T$ has two ${\cal H}$-classes and these are also
${\cal D}$-classes: $H_{\beta}=\{\beta^{r}:\,0\leq r\leq m-1\}$
of cardinal $m$, which is the group of units of $T$ and $H_{\alpha}=\{\beta^{r}\alpha\beta^{s}:0\leq r,s\leq m-1\}$
of cardinal $m^{2}$ and $H_{\alpha}<_{{\cal J}}H_{\beta}$. All members
$\gamma=\beta^{r}\alpha\beta^{s}$ of $H_{\alpha}$ are not regular;
dom $\gamma\subseteq\{X\times(m_{n+i}-r)\,\mbox{(mod \ensuremath{m)\,(0\leq i\leq n-1)\}}}$
with dom $\gamma$ meeting each specified interval and ran $\gamma\subseteq\{X\times(m_{i}+s)\,\mbox{(mod \ensuremath{m})\,\ensuremath{(0\leq i\leq n-1)}\ensuremath{\}}}$
with ran $\gamma$ similarly meeting each specified interval.

(ii) $T_{1}=\{\lambda(\alpha_{i},j,k):0\leq i\leq n-1,\,0\leq j,k\leq m-1\}\cup\{0\}$
is isomorphic to the Rees matrix semigroup $M={\cal M}^{0}[S,m,m,I_{m}]$,
where $I_{m}$ is the $m\times m$ identity matrix. Moreover $T_{1}$
is isomorphic to $(S\times B)/I$, where $B$ is the $m\times m$
combinatorial Brandt semigroup and $I$ is the ideal $S\times\{0\}$
of $S\times B$. For each $j\in\mathbb{Z}_{m},$ the set $T_{1,j}=\{\lambda(\alpha_{i},j,j):0\leq i\leq n-1\}$
is a subsemigroup of $T$ isomorphic to $S$. 

(iii) For any $\gamma\in T$, with dom $\gamma\cap(X\times\{j\})\neq\emptyset$,
$\gamma|_{X\times\{j\}}=\lambda_{i,j,k}$ for some $0\leq i\leq n-1,\,0\leq k\leq m-1$.

(iv) $T=T_{1}\cup H_{\alpha}\cup H_{\beta}$, and the union is a disjoint
union. Moreover\textbf{ $T_{1}$ }is an ideal of $T$ and if $S$
is regular then so is $T_{1}$.

(v) The set of idempotents $E(T)=\bigcup_{i=1}^{m}E_{i}\cup\{0,\iota\}$,
where $E_{i}=\{\lambda(e,i,i):e\in E(S),0\leq i\leq m-1\}$. Moreover
all products of non-identity idempotents equal $0$ except those within
some $E_{i}$. In particular if $E(S)$ is a band then so is $E(T)$. 

\emph{Proof }(i) The powers of $\beta$ are exactly the members of
$T$ with range (and domain) $Z$, and by Lemma 3.2 $\langle\beta\rangle$
is a cyclic group, the group of units of $T$, whence it follows that
$D_{\beta}=H_{\beta}=\langle\beta\rangle$ and by definition $|H_{\beta}|=m$.

The set $A=\{\beta^{r}\alpha\beta^{s}:r,s\geq0\}\subseteq H_{\alpha}$.
By Lemma 3.3, any $\delta\leq_{{\cal J}}\gamma$, where $\gamma\in T\setminus(A\cup H_{\beta})$
has domain within some single interval of $Z.$ If $\gamma\in D_{\alpha}$
we would have $\alpha\leq_{{\cal J}}\gamma$, whence dom $\alpha$
is contained in a single interval of $Z$, which contradicts the definition
of $\alpha$. It follows that $D_{\alpha}\subseteq A\subseteq H_{\alpha}\subseteq D_{\alpha}$,
giving equality throughout and $H_{\alpha}<_{{\cal J}}H_{\beta}$.

Next take $\gamma=\beta^{r}\alpha\beta^{s}$ so that 
\[
\mbox{dom \ensuremath{\gamma=\,\,}dom \ensuremath{\beta^{r}\alpha\beta^{s}=\{(x,(j-r)\,\mbox{(mod\,\ensuremath{m):\,(x,j)\in\,\mbox{dom\,\ensuremath{\alpha}\}. }}}}}
\]
Since dom $\alpha\subseteq\{(X,m_{n+i}):\,(0\leq i\leq n-1)\}$ and
dom $\alpha$ meets each of these intervals, it follows that dom $\gamma\subseteq\{(X,(m_{n+i}-r)\,(\mbox{mod \ensuremath{m})\,\ensuremath{(0\leq i\leq n-1)\}}}$
as stated and that dom $\gamma$ meets each of these intervals. Since
$\alpha$ maps the members of its domain in the interval $(X,m_{n+i})$
into the interval $(X,m_{i})$, the claim for ran $\gamma$ now follows
in the same way.

Suppose that $\gamma=\beta^{r_{1}}\alpha\beta^{s_{1}},\delta=\beta^{r_{2}}\alpha\beta^{s_{2}}$
and that $\gamma=\delta$. We wish to show that $\beta^{r_{1}}=\beta^{r_{2}}$
and $\beta^{s_{1}}=\beta^{s_{2}}$. By cancelling powers of $\beta$
in the equation of any counter example to this claim we would obtain
a counter example where $\gamma=\beta^{r}\alpha\beta^{s}$ and where
$\delta=\alpha$, $(0\leq r,s\leq m-1)$ so let us assume this case.
However since $|S|\geq3$ we have by Lemma 2.5 and our statement on
domains that dom $\gamma=$ dom $\alpha$ implies that $r=0$ and
similarly we have ran $\gamma=$ ran $\alpha$ implies $s=0$, as
required. We conclude that all products $\beta^{r}\alpha\beta^{s}$
$(0\leq r,s\leq m-1)$ are pairwise distinct and $|H_{\alpha}|=m^{2}$
as claimed. 

If any member of $D_{\alpha}$ were regular, the same would be true
of $\alpha$. However, by Lemmas 3.2 and 3.3, for any $\gamma\in T$
we have $\alpha\gamma\alpha\not\in D_{\alpha}$, so in particular
$\alpha=\alpha\gamma\alpha$ is impossible in $T$ and hence $D_{\alpha}$
is not a regular ${\cal D}$-class. 

(ii) From Lemma 3.5 and the definitions of $\alpha$ and $\beta$
we have the following formulae:
\begin{equation}
\lambda(\alpha_{i_{1}},j_{1},k)\lambda(\alpha_{i_{2}},k,k_{2})=\lambda(\alpha_{i_{1}}\alpha_{i_{2}},j_{1},k_{2})
\end{equation}
\begin{equation}
\lambda(\alpha_{i_{1}},j_{1},k_{1})\lambda(\alpha_{i_{2}},j_{2},k_{2})=0\,\,\mbox{if\,\,\ensuremath{k_{1}\neq j_{2}}}
\end{equation}
\begin{equation}
\beta\lambda(\alpha_{i},j,k)=\lambda(\alpha_{i},j-1,k),\,\,\lambda(\alpha_{i},j,k)\beta=\lambda(\alpha_{i},j,k+1)
\end{equation}
\begin{equation}
\alpha\lambda(\alpha_{i},m_{j},k)=\lambda(\alpha_{j}\alpha_{i},m_{n+j},k)\,(0\leq j\leq n-1)
\end{equation}
\begin{equation}
\alpha\lambda(\alpha_{i},j,k)=0\,\,\mbox{if\,\,\ensuremath{j\not\in\{m_{t}:0\leq t\leq n-1\}}}
\end{equation}
\begin{equation}
\lambda(\alpha_{i},j,m_{n+k})\alpha=\lambda(\alpha_{i}\alpha_{k},j,m_{k})
\end{equation}
\begin{equation}
\lambda(\alpha_{i},j,k)\alpha=0\,\,\mbox{if\,\,\ensuremath{k\not\in\{m_{n+t}:0\leq t\leq n-1\}}}
\end{equation}
From (4) and (5) we see that products in $T_{1}$ are indeed those
of the Rees matrix semigroup $M$, which is then isomorphic to $(S\times B)/I$.
The diagonal ${\cal H}$-classes of $M$ are each copies of our monoid
$S$. 

(iii) The claim is clearly true for $\gamma=\alpha,\beta$ as 
\[
\alpha|_{X\times\{m_{n+i}\}}=\lambda(\alpha_{i},m_{n+i},m_{i})\mbox{\,\ and \ensuremath{\beta|_{X\times\{i\}}=\lambda(\alpha_{0},i,i+1).}}
\]
The result now follows by induction on the length of $\gamma$ (taken
as a product in the generators $\alpha$ and $\beta$): let $\gamma=\rho\mu$
say, where $\mu\in\{\alpha,\beta\}$. Then $(\rho\mu)|_{(X\times\{j\})}=\rho|_{X\times\{j\}}\mu$
but by induction we may write this product as $\lambda(\alpha_{i},j,k_{1})\mu$
say. By formulae (6),(9), and (10) this in turn may be written as
$\lambda(\alpha_{l},j,k)|_{X\times\{j\}}=\lambda(\alpha_{l},j,k)$
for some $0\leq l\leq n-1$ and $0\leq k\leq m-1$, as required. 

(iv) Since the domains of members of $T_{1}$ are each contained within
a single interval and those of $H_{\alpha}\cup H_{\beta}$ are not,
we have by this and part (i) that the three sets are pairwise disjoint.
It remains to verify that if $\gamma\in T\setminus(H_{\alpha}\cup H_{\beta})$
then $\gamma\in T_{1}$. However, by Lemma 3.3 we have dom $\gamma\subseteq X\times\{j\}$
say and so by part (iii) we have either $\gamma=0$ or $\gamma=\gamma|_{X\times\{j\}}=\lambda_{i,j,k}$
for some $i,k$. In other words, $\gamma\in T_{1}$. From equations
(6 - 10) it follows that $T_{1}$ is an ideal of $T$. Finally for
any non-zero $\lambda=\lambda(\alpha_{i},j,k)\in T_{1}$ we have $\lambda(\alpha_{i}',k,j)\in T_{1}$
is an inverse of $\lambda$ in $T_{1}$ for any choice of $\alpha_{i}'\in V(\alpha_{i})$. 

(v) By (i), $\iota$ is the unique idempotent in $H_{\alpha}\cup H_{\beta}$.
Hence any other non-zero idempotent $\varepsilon$ belongs to $T_{1}$
and in particular dom $\varepsilon\subseteq X\times\{i\}$ say. Since
$\varepsilon$ is a non-zero idempotent, it follows that $\emptyset\neq Z\varepsilon\subseteq X\times\{i\}$.
Hence by (iv) we obtain $\varepsilon=\lambda(e,i,i)$ for some $e\in S$,
and clearly $e=e^{2}$ so that $\varepsilon\in E_{i}$, as claimed.
The claims regarding products of idempotents now follows. This completes
the proof of the theorem. $\qed$

\textbf{Corollary 3.7} Let $S$ be a finite monoid such that $E(S)$
is a subsemigroup of $S$. Then $S$ may be embedded in a finite monoid
$T$ such that $E(T)$ is a submonoid of $T$ and $T$ is generated
as a semigroup by a set of two generators $\{\alpha,\beta\}$ where
$\beta$ is a group element and $\alpha$ is nilpotent of index $2$.
Moreover if $|E(S)|\geq2$, then $E(T)$ satisfies the same semigroup
identities as $E(S)$. 

\textbf{Remark 3.8} If $|E(S)|\leq1$ then, since $S$ is a monoid
and every member of $S$ has an idempotent power, it follows that
$S$ is a finite group. We may then embed $S$ in the finite symmetric
group $T=G_{S}$, which is two-generated and then $E(S)$ and $E(T)$
are both trivial and so satisfy every semigroup identity.

\emph{Proof }Take $T=\langle\alpha,\beta\rangle$ as in Theorem 3.6.
It remains only to verify that if $\phi:p=q$ is a semigroup identity
satisfied by $E(S)$ then $\phi$ is satisfied by $E(T)$, the converse
implication being clear as $E(S)$ is embedded in $E(T)$. If one
side of $\phi$, the word $p$ say, had a variable $y$ that did not
appear in $q$, then substituting all other variables in $\phi$ by
$\iota$ gives the identity $y=1$, whence it follows that the monoid
$E(S)$ is trivial, contrary to hypothesis. Hence each variable $x$
of $\phi$ appears in both $p$ and $q$. 

By Theorem 3.6(v), all products of non-identity idempotents within
$E(T)$ equal $0$ unless they take place within some $E_{i}=\{\lambda(e,i,i):e\in E(S),0\leq i\leq m-1\}$.
Hence if, under some substitution from $E(T)$, one side of $\phi$,
$p$ say, is not $0$, then all variables of $\phi$ have been substituted
by either $\iota$ or by members of some subsemigroup $E_{i}$ of
$E(T)$. By replacing $\iota$ with the identity of $E_{i}$ as required,
we express the products $p$ and $q$ as products of members of $E_{i}$
while retaining the same values. However, since $E_{i}\cong E(S)$,
it follows that $p=q$ is satisfied in $E_{i}$ as well and so the
products $p$ and $q$ in $E_{i}$ are equal. It follows that $E(T)$
also satisfies the identity $\phi$. $\qed$

\textbf{Remark 3.9} In the case of a finite semigroup $S$ that is
not a monoid we may work with $S^{1}$. If $E(S)$ forms a band then
so does $E(S^{1})$ and the previous construction then yields a finite
2-generated monoid $T$ containing $S^{1}$ (and so containing $S$)
such that $E(T)$ is also a band.

\section{Orthodox semigroups }

We next use the construction of Section 3 to provide another proof
of Theorem 1.3 and to show that if the original semigroup $S$ is
orthodox, the same is true for the 2-generated containing semigroup
$T$. We will however now put $m_{i}=2^{i},\,i=0,1,2\cdots,2n-1$
so our modulus used for our cycle $\beta$ becomes $m=1+2^{2n-1}$.
Let $S=\{\alpha_{0},\alpha_{1},\cdots,\alpha_{n-1}\}$ now denote
a finite regular monoid with $\alpha_{0}=\iota$ and $S\leq PT_{X}$
for some finite base set $X$ as before. We may also assume that the
domain of each mapping $\alpha_{i}$ is not empty.

For each $\alpha_{i}\in S$ choose and fix an inverse $\alpha_{i}'\in V(\alpha_{i})$
(there is no assumption that the mapping $(')$ on $S$ is one-to-one).
The cycle $\beta$ is just as before and its action is given by (1).
Similarly, the action (2) remains applicable to our second generator
$\alpha$. However we augment the domain of $\alpha$ to include all
the intervals $X\times\{m_{i}\}$ $(0\leq i\leq n-1)$, the union
of which contained the range set of $\alpha$ but previously lay outside
of the domain of $\alpha$. Define:
\begin{equation}
(x,m_{i})\cdot\alpha=(x\cdot\alpha_{i}',m_{n+i})\,\,(0\leq i\leq n-1)
\end{equation}

\textbf{Remarks 4.2 }It will be convenient to also denote $\alpha_{i}'$
by $\alpha_{i+n}$, in which case the definition of the action of
$\alpha$ is encapsulated by:
\begin{equation}
(x,m_{t})\cdot\alpha=(x\cdot\alpha_{t\pm n},m_{t\pm n})\,\,(0\leq t\leq2n-1)
\end{equation}
where the signs associated with the $\pm$ signs in (12) are not independent
but are equal to each other: the sign on the subscripts is $+$ or
$-$ according as $0\leq t\leq n-1$ or $n\leq t\leq2n-1$. Although
$\alpha$ is no longer a nilpotent (see Lemma 4.3) it is still the
case that any $\gamma\in T$ acts in a one-to-one fashion on the second
entries of the pairs $(x,i)\in$ dom $\gamma$ (as shown in the proof
of Lemma 3.2) and $\gamma$ maps intervals into intervals as this
holds for each of the generators $\alpha$ and $\beta$. We next prove
the counterpart of Lemma 3.3. 

\textbf{Lemma 4.3}

(i) The mappings $\alpha$ and $\beta$ of $T$ satisfy $\beta=\beta^{m}$
and $\alpha=\alpha^{3}$.

(ii) Let $\gamma=\beta^{r}\alpha^{\varepsilon}\beta^{s}$ for $\varepsilon\geq1$.
Then \newline dom $\gamma$ $\subseteq\{X\times\{(m_{t}-r)\,\mbox{(mod \ensuremath{m),\,(0\leq t\leq2n-1)\}}}$
and dom $\gamma$ has non-empty intersection with each of these intervals.
Similarly ran $\gamma\subseteq\{X\times(m_{t}+s)\mbox{ (mod \ensuremath{m)\,(0\leq t\leq2n-1)\}}}$
with ran $\gamma$ meeting each of these intervals.

(iii) Let $\gamma\leq_{{\cal J}}\alpha\beta^{r}\alpha$ where $r\not\equiv0$
(mod $m$). Then dom $\gamma\subseteq X\times\{i\}$ for some $0\leq i\leq m-1$.

\emph{Proof }(i) That $\beta=\beta^{m}$ is true as before. For any
$(x,m_{n+i})\in$ dom $\alpha$ we have by (12) that 
\[
(x,m_{n+i})\cdot\alpha^{3}=(x\cdot\alpha_{i},m_{i})\cdot\alpha^{2}=(x\cdot\alpha_{i}\alpha_{i}',m_{n+i})\cdot\alpha
\]
\[
=(x\cdot\alpha_{i}\alpha_{i}'\alpha_{i},m_{i})=(x\cdot\alpha_{i},m_{i})=(x,m_{n+i})\cdot\alpha,
\]
and in the same way we obtain $(x,m_{i})\cdot\alpha^{3}=(x,m_{i})\cdot\alpha,$
thus showing that $\alpha=\alpha^{3}$. Note also that by finiteness
it follows that $\alpha|_{\mbox{ran \ensuremath{\alpha}\ }}$ is a
permutation and so dom $\alpha=$ dom $\alpha^{2}$ and ran $\alpha=$
ran $\alpha^{2}$. 

(ii) Let us write (for the purposes of this part only) 
\[
D_{\gamma}=\{i:(X\times\{i\})\cap\,\mbox{dom\,\ensuremath{\gamma\neq\emptyset\}\,}and \ensuremath{R_{\gamma}=\{i:(X\times\{i\})\cap\,\mbox{ran \ensuremath{\gamma\neq\emptyset\}.}}}}
\]
Observe that for any $\varepsilon\geq1$, $\ensuremath{D_{\alpha^{\varepsilon}}=R_{\alpha^{\varepsilon}}=\{m_{t}:0\leq t\leq2n-1\}}$.
Also note that for any $\gamma\in T$ we have $D_{\beta^{r}\gamma\beta^{s}}=(D_{\gamma}-r)$
(mod $m$) and $R_{\beta^{r}\gamma\beta^{s}}=(R_{\gamma}+s)$ (mod
$m$). Applying these facts to $\gamma=\alpha^{\varepsilon}$ then
proves the claims of (ii). 

(iii) As in the proof of Lemma 3.3, it is enough to consider the case
represented by $\gamma=\alpha\beta^{r}\alpha$. Since dom $\gamma\subseteq$
dom $\alpha$, it follows that each member of dom $\gamma$ has the
form $(x,m_{t})$ for some $0\leq t\leq2n-1$ and so 
\begin{equation}
(x,m_{t})\cdot\alpha\beta^{r}\alpha=(x\cdot\alpha_{t\pm n},m_{t\pm n})\cdot\beta^{r}\alpha=(x\cdot\alpha_{t\pm n},m_{t\pm n}+r\,\mbox{(mod \ensuremath{m)}})\cdot\alpha
\end{equation}
This implies that $m_{t\pm n}+r\equiv m_{k}$ (mod $m$) for some
$0\leq k\leq2n-1$. Now suppose that $(x',m_{t'})\in$ dom $\gamma$;
by (13) we deduce that $m_{t'\pm n}+r\equiv m_{k'}$ (mod $m$) for
some $0\leq k'\leq2n-1$, which yields: 
\begin{equation}
m_{t\pm n}-m_{k}\equiv m_{t'\pm n}-m_{k'}\equiv-r\,\mbox{(mod \ensuremath{m}) }
\end{equation}
where the signs taken in the $\pm$ symbols occurring in (14) are
not necessarily equal to each other. If the first congruence in (14)
is equality then since $r\not\equiv0$ (mod $m$), we have that $m_{t\pm n}\neq m_{k}$
and $m_{t'\pm n}\neq m_{k'}$ and so by the MC property $m_{t\pm n}=m_{t'\pm n}$
(and $m_{k}=m_{k'}$). It follows either that $t=t'$ or 
\[
((t-n=t'+n)\,\mbox{or\,\ensuremath{(t+n=t'-n))\Rightarrow|t-t'|=2n.}}
\]
However, since $0\leq t,t'\leq2n-1$, the latter is not possible and
so $t=t'$. Otherwise the congruence in (14) is not equality whence:
\begin{equation}
(m_{t\pm n}+m_{k'})-(m_{t'\pm n}+m_{k})=\pm(1+m_{2n-1})
\end{equation}
By multiplying throughout by $-1$ and interchanging $t$ and $t'$
if necessary, we may take the $+$ sign in (15). Since $r\not\equiv0$
(mod $m$) we have that $m_{t\pm n}\neq m_{k}$ and $m_{t'\pm n}\neq m_{k'}$.
However, by Lemma 2.3, one term in the first bracket equals $m_{2n-1}$,
one term in the second bracket equals $1$ and the other two terms
cancel each other. 

Hence either $m_{t\pm n}=2^{2n-1},$ $m_{t'\pm n}=1$ and $m_{k}=m_{k}'$,
or $m_{k}'=2^{2n-1},$ $m_{k}=1$ and $m_{t\pm n}=m_{t'\pm n}$. However
$m_{k}=m_{k'}$ implies (by (14)) that $m_{t\pm n}=m_{t'\pm n}$ and
so $t=t'$ is the conclusion. Similarly the latter possibility once
again gives $t=t'$. Therefore dom $\gamma\subseteq X\times\{m_{t}\}.$
$\qed$ 

Lemmas 3.4 and 3.5 are valid for our extended construction, the proofs
being unchanged from the originals. Moreover the description of the
mapping $\lambda_{i,j,k}$ of Theorem 3.6(ii) continues to hold in
our monoid $T$ currently under consideration, as do the formulae
(4 - 6). The full set of corresponding formulae for $T$ (additions
and subtractions taken mod $m$) are as follows:

\begin{equation}
\beta\lambda(\alpha_{i},j,k)=\lambda(\alpha_{i},j-1,k),\,\,\lambda(\alpha_{i},j,k)\beta=\lambda(\alpha_{i},j,k+1)
\end{equation}
\begin{equation}
\alpha\lambda(\alpha_{i},m_{j},k)=\lambda(\alpha_{j}\alpha_{i},m_{j\pm n},k)\,(+\,\mbox{if \ensuremath{0\leq j\leq n-1,\,-\mbox{\,\ if\,\ \ensuremath{n\leq j\leq2n-1)}}}}
\end{equation}
\begin{equation}
\alpha\lambda(\alpha_{i},j,k)=0\,\,\mbox{if\,\,\ensuremath{j\not\in\{m_{t}:0\leq t\leq2n-1\}}}
\end{equation}
\begin{equation}
\lambda(\alpha_{i},k,m_{j})\alpha=\lambda(\alpha_{i}\alpha_{j\pm n},k,m_{j\pm n})\,(+\,\mbox{if \ensuremath{0\leq j\leq n-1,\,-\,\mbox{if \ensuremath{n\leq j\leq2n-1}}}})
\end{equation}
\begin{equation}
\lambda(\alpha_{i},j,k)\alpha=0\,\,\mbox{if\,\,\ensuremath{k\not\in\{m_{t}:0\leq t\leq2n-1\}}}
\end{equation}

\textbf{Proposition 4.4 }Let $T=\langle\alpha,\beta\rangle$.

(i) For any $\gamma\in T$, with dom $\gamma\cap(X\times\{j\})\neq\emptyset$,
$\gamma|_{X\times\{j\}}=\lambda_{i,j,k}$ for some $0\leq i\leq n-1,\,0\leq k\leq m-1$;

(ii) $T$ is regular. 

\emph{Proof }(i) The claim is clearly true for $\gamma=\alpha,\beta$
as 
\[
\alpha|_{X\times\{m_{i}\}}=\lambda(\alpha_{i\pm n},m_{i},m_{i\pm n})\,\mbox{and\,\ensuremath{\beta|_{X\times\{i\}}=\lambda(\alpha_{0},i,i+1).}}
\]
The result now follows as in Theorem 3.6 (iii) by induction on the
length of $\gamma$ (taken as a product in the generators $\alpha$
and $\beta$), together with formulae (16 - 20).

(ii) Take an arbitrary product $p=\beta^{r_{1}}\alpha\beta^{r_{2}}\alpha\cdots\beta^{r_{t-1}}\alpha\beta^{r_{t}}\in T$
with $(1\leq t,\,0\leq r_{i}\leq m-1)$. If $t=1$, then $p=\beta^{r_{1}}$
is a group element and so $p$ is regular. Since $\alpha=\alpha^{3}$
it follows that all mappings of the form $\beta^{r}\alpha^{\varepsilon}\beta^{s}$
$(\varepsilon=1,2)$ are contained in the regular ${\cal D}$-class
$D_{\alpha}$ of $T$. This deals with the case where $t=2$ and the
case $(t=3$ and $r_{2}=0$). The remaining cases are where $t\geq3$
and $p$ has one of the two forms $p=\beta^{r_{1}}\alpha\beta^{r_{2}}\alpha\beta^{r_{3}}\cdots$
or $p=\beta^{r_{1}}\alpha^{2}\beta^{r_{2}}\alpha\beta^{r_{3}}\cdots$
with $r_{2}\not\neq0$ in both instances. It follows from Lemma 4.3(iii)
that dom $p\subseteq X\times\{j\}$ say. Of course if $p=0$ then
$p$ is regular. Otherwise by (i) $p=p|_{X\times\{j\}}=\lambda_{i,j,k}$
for some $0\leq i\leq n-1,$ $0\leq k\leq m-1$. By Theorem 3.6(ii),
$p$ is a member of a subsemigroup of $T$ isomorphic to $(S\times B)/I$,
and in particular $p$ is a regular member of $T$. $\qed$

Proposition 4.4 shows that any finite semigroup may be embedded in
a finite regular semigroup $T$ generated by two group elements, thereby
providing a new proof of Theorem 1.3. However, the semigroup $T$
preserves the idempotent structure of $S$ in that $E(T)$ consists
of copies of $E(S)$ together with the conjugates under $\beta$ of
$\alpha^{2}$. 

\textbf{Theorem 4.5 }(Structure of $T$) 

(i) $H_{\beta}$ is the group of units of $T$, which is cyclic of
order $m$. Moreover $D_{\alpha}<_{{\cal J}}H_{\beta}$ and $D_{\alpha}=\{\beta^{r}\alpha^{\varepsilon}\beta^{s}:\,\varepsilon=1,2\}$. 

(ii) The monoid $T$ has an ideal $T_{1}$ with $\gamma<_{{\cal J}}\alpha$
for all $\gamma\in T_{1}$ where $T_{1}=\{\lambda_{i,j,k}\}\cup\{0\}$
$(0\leq i\leq n-1,0\leq j,k\leq m-1)$. 

(iii) $T=H_{\beta}\cup D_{\alpha}\cup T_{1}$ with the union a disjoint
union.

(iv) The set of idempotents of $T$ is given by $E(T)=E\cup F\cup\{\iota,0\}$,
where $E=\{\lambda(e,i,i):e\in E(S),\,0\leq i\leq m-1\}$ and $F=\{\beta^{j}\alpha^{2}\beta^{-j}:0\leq j\leq m-1\}$.
Moreover each $\rho\in E(T)$ maps identically on its second entry,
meaning that $(X\times\{i\})\rho\subseteq X\times\{i\}$. 

(v) The principal factor $D_{\alpha}\cup\{0\}$ of $T$ is of cardinal
$1+2m^{2}$ and is a Brandt semigroup ${\cal M}^{0}[\mathbb{Z}_{2},m,m,I_{m}]$. 

\emph{Proof }(i) As in Section 3, $H_{\beta}$ is the group of units
of $T$ of cardinal $m$. Also $\gamma<_{{\cal J}}\beta$ for any
$\gamma\in S\alpha S$ and so $D_{\alpha}<_{{\cal J}}H_{\beta}$.
By Lemma 4.3(i), $\alpha=\alpha^{3}$ and so $A=\{\beta^{r}\alpha^{\varepsilon}\beta^{s}:\,\varepsilon=1,2\}\subseteq D_{\alpha}$.
Conversely, if $\alpha\leq_{{\cal J}}\gamma$ with $\gamma\in T\setminus(A\cup H_{\beta})$
then $\alpha\leq_{{\cal J}}\gamma\leq_{{\cal J}}\alpha\beta^{r}\alpha$
for some $r\not\equiv0$ (mod $m$) and by Lemma 4.3(iii), it would
follow that dom $\alpha$ was contained in a single interval of $T$,
contrary to the definition of $\alpha$. Hence $A=D_{\alpha}$, thus
establishing (i).

(ii) As in the proof of Lemma 3.5, we have that $T_{1}\subseteq T$
and that $T_{1}$ is an ideal of $T$ follows from the formulae (16
- 20). From Lemma 3.5 we have that $\gamma\not\in H_{\beta}$ whence
$\gamma\leq_{{\cal J}}\alpha$ and that the inequality is strict follows
from Proposition 4.4(i) and the fact that, unlike dom$\gamma$, dom
$\alpha$ is not contained in a single interval.

(iii) It follows from parts (i) and (ii) that $H_{\beta}\cup D_{\alpha}\cup T_{1}\subseteq T$
and the union is a disjoint union. Conversely take any $\gamma\in T\setminus\{H_{\beta}\cup D_{\alpha}\}$.
By part (i), Lemma 4.3(iii) applies to $\gamma$ whence by Proposition
4.4(i) it follows that $\gamma\in T_{1}$, as required. 

(iv) Clearly all the members listed in $E(T)$ are indeed idempotents.
For any $\lambda=\lambda_{i,j,k}\in T_{1}$ we have $\lambda^{2}=0$
unless $k=j$, in which case $\lambda^{2}=\lambda$ if and only if
$\alpha_{i}=e\in E(S)$ and so $\lambda=\lambda(e,j,j)\in E$. From
part (iii) it follows that all other members $p\in E(T)$, other than
$0$ and $\iota$, lie in $D_{\alpha}$ and so have the form $p=\beta^{j}\alpha^{\varepsilon}\beta^{k}$
where $(\varepsilon\in\{1,2\})$. We next check that if $j+k\equiv0$
(mod $m$) then $p=p^{2}$ if and only if $\varepsilon=2$. The reverse
implication just says that all members of $F$ are idempotents, which
has already been noted, so let us suppose that, contrary to our claim,
$\varepsilon=1$ and we have $p=\beta^{j}\alpha\beta^{-j}$ with that
$p=p^{2}$. Then $\beta^{j}\alpha\beta^{-j}=\beta^{j}\alpha^{2}\beta^{-j}$,
which in turn implies that $\alpha=\alpha^{2}$, which is false as
$X\times\{m_{n}\}$ is an interval that meets dom $\alpha=$ dom $\alpha^{2}$
but $(X\times\{m_{n}\})\alpha\subseteq X\times\{m_{0}\},\,(X\times\{m_{n}\})\alpha^{2}\subseteq X\times\{m_{n}\}$. 

Let us therefore examine the case where $j+k\not\equiv0$ (mod $m$)
for some $0\leq j,k\leq m-1$. Since $p=p^{2}$ and the product $p^{2}$
contains a factor of the form $\alpha\beta^{t}\alpha$ with $t\not\equiv0$
(mod $m$), it now follows by Lemma 4.3(iii) and the fact that $p=p^{2}$
that both dom $p$ and ran $p$ are contained in $X\times\{i\}$ say.
However, since dom $\alpha=$ dom $\alpha^{2}=X\times\{m_{t}\}_{0\leq t\leq2n-1}$,
it follows from Lemma 4.3(ii) that dom $p=$ dom $\beta^{j}\alpha^{\varepsilon}\beta^{k}=X\times\{m_{t}-j\}_{0\leq t\leq2n-1}$
$(\varepsilon\in\{1,2\})$. In particular, dom $p$ is not contained
within a single set of the form $X\times\{i\}$ and this contradicts
the assumption that $p\in E(T)$. Therefore the set $E(T)$ is as
described. The final assertion is clearly true for idempotents $0$
and $\iota$ and those in $T_{1}$. By above, any idempotent $\rho\in F$
satisfies dom $\rho\subseteq X\times\{i\}$ say and since any idempotent
maps identically on its range it follows that $(X\times\{i\})\rho\subseteq X\times\{i\}$
from which the claim follows.

(v) There are $2m^{2}$ expressions of the form $\beta^{r}\alpha^{\varepsilon}\beta^{s}:(\varepsilon\in\{1,2\},\,0\leq r,s\leq m-1)$
and so the cardinality claim will follow by showing they are pairwise
distinct. If not, we would have an equality of the form $\alpha^{\varepsilon_{1}}=\beta^{r}\alpha^{\varepsilon_{2}}\beta^{s}=\gamma$
say, for some $\varepsilon_{1},\varepsilon_{2}\in\{1,2\}$. By Lemma
4.3(ii), dom $\gamma$ $\subseteq\{X\times\{(m_{t}-r)\,\mbox{(mod \ensuremath{m),\,0\leq t\leq2n-1\}}}$
and dom $\gamma$ has non-empty intersection with each of these intervals.
Since $|S|\geq3$ it follows by Lemma 2.5 that $r=0$ and in the same
way we infer likewise that $s=0$ as well.

Since $D_{\alpha}$ is a regular ${\cal D}$-class, the principal
factor $P=D_{\alpha}\cup\{0\}$ is a completely 0-simple semigroup.
By part (i) and Lemma 4.3 parts (ii) and (iii) we see that for $\gamma=\beta^{r}\alpha^{\varepsilon}\beta^{s}\in D_{\alpha}$
we have $R_{\gamma}=\{\beta^{r}\alpha^{\varepsilon}\beta^{t};\varepsilon\in\{1,2\},\,0\leq t\leq m-1\}$,
$L_{\gamma}=\{\beta^{t}\alpha^{\varepsilon}\beta^{s},\varepsilon\in\{1,2\},\,0\leq t\leq m-1\}$
and so $H_{\gamma}=\{\beta^{r}\alpha^{\varepsilon}\beta^{s}:\varepsilon\in\{1,2\}\}.$
In particular $H_{\alpha}=\{\alpha,\alpha^{2}\}\cong\mathbb{Z}_{2}$.
By the previous paragraph it follows that there are $m$ ${\cal R}$-classes
and $m$ ${\cal L}$-classes of $D_{\alpha}$, so that $P\cong$${\cal M}^{0}[\mathbb{Z}_{2},m,m,M]$
is the Rees matrix form of this principal factor for some $m\times m$
matrix $M$. To complete the proof we only need to know that the idempotents
of $P$ form a semilattice, for then $P$ is a regular $0$-simple
semigroup with commuting idempotents, which is necessarily a Brandt
semigroup, whence $M$ can be taken to be the identity matrix. However,
the product of any two distinct idempotents $e=\beta^{j}\alpha^{2}\beta^{-j}$
and $f=\beta^{k}\alpha^{2}\beta^{-k}$ is $\beta^{j}\alpha^{2}\beta^{k-j}\alpha^{2}\beta^{k}$
and since $k\not\equiv j$ (mod $m$) it follows from (i) above together
with Lemma 4.3(iii) that $ef\not\in D_{\alpha}$ so that \emph{in
the principal factor} $D_{\alpha}\cup\{0\}$, the product of any two
distinct idempotents is $0$ and in particular $E(D_{\alpha}\cup\{0\})$
is a semilattice, as required. $\qed$

\textbf{Theorem 4.6} (a) Any finite orthodox semigroup $S$ may be
embedded in a finite orthodox semigroup $T$ generated by two group
elements.

(b) Any finite orthodox monoid $S^{1}$ may be embedded as a semigroup
into a finite 2-generated orthodox monoid $T$ whose subband of idempotents
satisfy the same semigroup identities.

\emph{Proof }(a) From Proposition 4.4, we need only check that, given
that $S$ is orthodox, the idempotents of our containing semigroup
$T$ form a band. Consider $E(T)=E\cup F\cup\{0\}$ as described in
Theorem 4.5. Products involving $0$ are $0$ and the product of any
two members of $E$ is also $0$ unless they have identical second
and third co-ordinates $j$ say. In this case we have a product of
idempotents in the semigroup $T_{1,j}\cong S$ by Theorem 3.6(ii):
in particular the product is itself an idempotent as $S$ is orthodox.

Next, let $\rho=\beta^{j}\alpha^{2}\beta^{-j}$ and $\mu=\beta^{k}\alpha^{2}\beta^{-k}$
be two distinct members of $F$. Since the product $\rho\mu$ has
the factor $\alpha\beta^{-j+k}\alpha$ with $k-j\not\equiv0$ (mod
$m$), it follows from Theorem 4.5 (iii) and (v) that either $\rho\mu=0$
or dom $\rho\mu\subseteq X\times\{i\}$ say. Routine calculation then
gives that, if defined, $(x,i)\cdot\rho\mu=(x\cdot ef,i)$ for some
idempotents $e,f\in E(S)$. Since $ef\in E(S)$ it follows that $\rho\mu=\lambda(ef,i,i)\in E(T)$.
In detail we have, working modulo $m$ with $i+j\equiv m_{t}$ (mod
$m$) say:
\[
(x,i)\cdot\rho=(x,i)\cdot\beta^{j}\alpha^{2}\beta^{-j}=(x,m_{t})\cdot\alpha^{2}\beta^{-j}=(x\cdot\alpha_{t\pm n},m_{t\pm n})\cdot\alpha\beta^{-j}
\]
\[
=(x\cdot\alpha_{t\pm n}\alpha_{t},m_{t})\beta^{-j}=(x\cdot\alpha_{t\pm n}\alpha_{t},i);
\]
now $\alpha_{t}$ is inverse to $\alpha_{t\pm n}$, so this final
product can be written as $(x$$\cdot e,i)$, where $e=\alpha_{t\pm n}\alpha_{t}\in E(S)$.
By the same token, applying this calculation now to $(x\cdot e,i)\cdot\mu$
yields the required expression $(x\cdot ef,i)$ where $ef\in E(S)$
as claimed previously. Hence $\rho\mu=\lambda(ef,i,i)\in E$. 

Finally let $\lambda=\lambda(e,i,i)\in E$ and $\rho=\beta^{j}\alpha^{2}\beta^{-j}\in F$
as above. If $\lambda\rho\neq0$ then $\lambda\rho$ has the form
$\lambda\rho=(ef,i,i)\in E(T)$ as $E(S)$ is a band. On the other
hand $\rho\lambda\neq0$ implies that $(x,i)\cdot\rho\lambda=(x\cdot fe,i,i)$
for some $f\in E(S)$ whence $\rho\lambda\in E(T)$. In detail the
relevant calculations are as follows. If $\lambda\rho\neq0$ then
dom $\lambda\rho\subseteq X\times\{i\}$, $i+j\equiv m_{t}$ (mod
$m$) say and 
\[
(x,i)\cdot\lambda\rho=(x\cdot e,i)\cdot\beta^{j}\alpha^{2}\beta^{-j}=(x\cdot e,m_{t})\cdot\alpha^{2}\beta^{-j}=(x\cdot e\alpha_{t\pm n},m_{t\pm n})\cdot\alpha\beta^{-j}=
\]
\[
(x\cdot e\alpha_{t\pm n}\alpha_{t},m_{t})\cdot\beta^{-j}=(x\cdot e\alpha_{t\pm n}\alpha_{t},i)
\]
and, as before, $\alpha_{t\pm n}\alpha_{t}=f\in E(S)$ and so $ef\in E(S)$
as $S$ is orthodox. Hence $(x,i)\cdot\lambda\rho=(x\cdot ef,i)$
and it follows that $\lambda\rho=\lambda(ef,i,i)\in E$. 

Now consider $\rho\lambda$ and suppose that $\rho\lambda\neq0$.
We have by Lemma 3.2 applied to $\rho$ that dom $\rho\lambda\subseteq X\times\{k\}$
say. However $(X\times\{k\})\rho$ meets dom $\lambda\subseteq X\times\{i\}$
and since $\rho$ is idempotent we have that $\rho$ maps each interval
$X\times\{k\}$ into itself and we deduce that $k=i$. Now we have
$i+j\equiv m_{t}$(mod $m$) say and we obtain: 
\[
(x,i)\cdot\rho\lambda=(x,i)\cdot\beta^{j}\alpha^{2}\beta^{-j}\lambda=(x,m_{t})\cdot\alpha^{2}\beta^{-j}\lambda=(x\cdot\alpha_{t\pm n},m_{t\pm n})\cdot\alpha\beta^{-j}\lambda=
\]
\[
(x\cdot\alpha_{t\pm n}\alpha_{t},m_{t})\beta^{-j}\lambda=(x\cdot\alpha_{t\pm n}\alpha_{t},i)\cdot\lambda=(x\cdot fe,i)
\]
 where $f=\alpha_{t\pm n}\alpha_{t}\in E(S)$ as before and again
$fe\in E(S)$ as $S$ is orthodox. Therefore $\rho\lambda=\lambda(fe,i,i)\in E$,
as required to complete the proof.

(b) Following Remark 3.8, only the case where $|E(S^{1})|\geq2$ is
of interest. As in the proof of Corollary 3.7, we may take a typical
semigroup identity $\phi:p=q$ satisfied by $S^{1}$ to be homotypical,
meaning that each variable in $\phi$ appears in both $p$ and $q$.
Since we are considering identities on bands, we may assume that $\phi$
has more than one variable. We need to check is that $E(T)$ also
satisfies $\phi$. 

By Lemma 4.3(iii) it follows that any product $uv$ of two distinct
members $u,v\in F=E(T)\cap D_{\alpha}$ falls out of $D_{\alpha}$
and lies in $T_{1}$. It follows, again from Lemma 4.3(iii), that
either $uv=0$ (the empty map) or dom$(uv)$, ran$(uv)$ are contained
in some interval $Y_{i}$ say. In the latter case $uv=(u|Y_{i})(v|Y_{i})$.
Since the restrictions $u_{i}=u|Y_{i}$ and $v_{i}=v|Y_{i}$ each
belong to $E_{i}=\{\lambda(e,i,i):e\in S\}$, the product $uv=u_{i}v_{i}$
is equal to a product of two idempotents in $E_{i}$.

Now let us consider the words $p(x_{1},\cdots,x_{t})$ and $q(x_{1},\cdots,x_{r})$$(r\geq2)$
of the identity $\phi$ and let us substitute elements of $E(T)$
to obtain products $P=p(t_{1},\cdots,t_{r})$ and $Q=q(t_{1},\cdots,t_{r})$.
We need to verify that $P=Q$. Since each product involves at least
$2$ members of $E(T)$, it follows from the argument of the previous
paragraph that each $t_{j}$ may be replaced by a member of $E(T_{1})$
without changing the value of either of the products $P$ and $Q$,
so without loss we may assume that $t_{1},\cdots,t_{r}\in E(T_{1})$.
Hence each $t_{j}\in E_{i}$ for some $i$ that depends on $j$. Consider
the set of subscripts $I=\{i:t_{j}\in E_{i}\}$. If $|I|=1$ then
both $P$ and $Q$ are products of idempotents in some $E_{i}\cong E(S)$
and so $P=Q$ as $E(S)$ satisfies $\phi$. On the other hand, if
$|I|\geq2$ then $P=Q=0$ as each of $P$ and $Q$ contains a product
of the form $uv$ with $u\in E_{i}$, $v\in E_{j}$ with $i\neq j$.
In either event, it follows that $\phi$ is satisfied by $E(T)$ also,
thus completing the proof of Theorem 4.6(b). $\qed$

Specialising to the case where $E(S)$ is a semilattice and noting
that $E(S)$ is a semilattice if and only if the same is true of $E(S^{1})$
gives the main corollary (Corollary 2.2) of the construction of {[}8{]}
that the finite symmetric inverse semigroup $I_{n}$ embeds in a $2$-generator
inverse susbsemigroup of $I_{n+2}$.

\textbf{Corollary 4.7 }(McAlister, Stephen and Vernitski) Every finite
inverse semigroup may be embedded in a finite 2-generated semigroup
that is an inverse semigroup.

\end{document}